\documentclass[11pt]{article}
\usepackage{amsfonts,latexsym,amsmath,amscd}
\usepackage[pdftex]{graphicx}
\newcommand{\cleqn}{\setcounter{equation}{0}}
\newcommand{\clth}{\setcounter{theorem}{0}}
\newcommand {\sectionnew}[1]{\section{#1}\cleqn\clth}
\newcommand{\beq}{\begin{equation}}
\newcommand{\eeq}{\end{equation}}
\newcommand{\beqa}{\begin{eqnarray}}
\newcommand{\eeqa}{\end{eqnarray}}
\newcommand{\beaa}{\begin{eqnarray*}}
\newcommand{\ben}{\begin{eqnarray*}}
\newcommand{\eaa}{\end{eqnarray*}}
\newcommand{\een}{\end{eqnarray*}}

\newcommand \nc {\newcommand}
\nc \proof {\noindent {\em{Proof.\/ }}} \nc \qed {$\Box$\hfill}
\newtheorem{theorem}{Theorem}[section]
\newtheorem{lemma}[theorem]{Lemma}
\newtheorem{proposition}[theorem]{Proposition}
\newtheorem{corollary}[theorem]{Corollary}
\newtheorem{definition}[theorem]{Definition}
\newtheorem{example}[theorem]{Example}
\newtheorem{remark}[theorem]{Remark}
\newtheorem{conjecture}[theorem]{Conjecture}
\newtheorem{question}[theorem]{Question}
\nc \bth[1] { \begin{theorem}\label{t#1} } \nc \ble[1] {
\begin{lemma}\label{l#1} } \nc \bpr[1] {
\begin{proposition}\label{p#1} } \nc \bco[1] {
\begin{corollary}\label{c#1} } \nc \bde[1] {
\begin{definition}\label{d#1}\rm } \nc \bex[1] {
\begin{example}\label{e#1}\rm } \nc \bre[1] {
\begin{remark}\label{r#1}\rm } \nc \bcon[1] {
\begin{conjecture}\label{con#1}\rm } \nc \bque[1] {
\begin{question}\label{que#1}\rm }
\nc {\eth} { \end{theorem} } \nc {\ele} { \end{lemma} } \nc {\epr} {
\end{proposition} } \nc {\eco} { \end{corollary} } \nc {\ede} {
\end{definition} } \nc {\eex} { \end{example} } \nc {\ere} {
\end{remark} } \nc {\econ} { \end{conjecture} } \nc {\eque} {
\end{question} }
\nc \thref[1]{Theorem \ref{t#1}} \nc \leref[1]{Lemma \ref{l#1}} \nc
\prref[1]{Proposition \ref{p#1}} \nc \coref[1]{Corollary \ref{c#1}}
\nc \deref[1]{Definition \ref{d#1}} \nc \exref[1]{Example \ref{e#1}}
\nc \reref[1]{Remark \ref{r#1}}

\def \R {{\mathcal R}}

\def \Rset {{\mathbb R}}

\def \Gr { {\mathrm{Gr}} }

\nc \GRA[1] { \Gr_A^{(#1)} }   
\nc \GRAN { \GRA{N} } \nc \GrA[1] { \Gr_A(#1) }\nc \GrAa {
\GrA{\alpha} }
\nc \GRB[1] { \Gr_B^{(#1)} }   
\nc \GRBN { \GRB{N} } \nc \GrB[1] { \Gr_B(#1) } \nc \GrBb {
\GrB{\beta} }
\nc \GRMB[1] { \Gr_{MB}^{(#1)} }   
\nc \GRMBN { \GRMB{N} } \nc \GrMB[1] { \Gr_{MB}(#1) } \nc \GrMBb {
\GrMB{\beta} }

\begin{document}
\title{{\Large\bf{Uniqueness of values of Aronsson operators and running costs in ``tug-of-war" games}}}
\author{Yifeng Yu\\
        Department of math, University of California at Irvine}
\date{}
\maketitle
\begin{abstract}
Let $A_H$ be the Aronsson operator associated with a Hamiltonian
$H(x,z,p).$ Aronsson operators arise from $L^\infty$ variational
problems, two person game theory, control problems, etc. In this
paper, we prove, under suitable conditions,  that if $u\in
W^{1,\infty}_{\rm loc}(\Omega)$ is  simultaneously a viscosity
solution of both of the equations \beq \label{a}
A_H(u)=f(x)\text{\quad and\quad } A_H(u)=g(x) \text{\quad in\quad}
\Omega, \eeq where $f, g\in C(\Omega),$ then $f=g.$ The assumption
$u\in W_{loc}^{1,\infty}(\Omega)$ can be relaxed to
 $u\in C(\Omega)$ in many interesting situations. Also,  we prove that if $f,g,u\in C(\Omega)$ and $u$
 is simultaneously a viscosity solution of the
equations \beq \label{b} {\Delta_{\infty}u\over
|Du|^2}=-f(x)\text{\quad and\quad} {\Delta_{\infty}u\over
|Du|^2}=-g(x) \text{\quad in\quad} \Omega,\eeq then $f=g.$ This
answers a question posed in Peres, Schramm, Scheffield and Wilson
[PSSW] concerning whether or not the value function uniquely
determines the running cost in the ``tug-of-war" game.
\end{abstract}
\sectionnew{Introduction}
Let $\Omega$ be an open set in $\Rset ^n$ and $H(x,z,p)\in
C^1(\Omega\times \Rset \times \Rset ^n)$. The Aronsson operator
associated with $H$ has the form
$$
A_{H}(u)=H_{p}(x,u,Du)\cdot D_{x}(H(x,u,Du)),
$$
where $D_{x}$ represents the partial derivative with respect to $x$
if we consider $H(x,u(x),Du(x))$ as a function of $x$. This type of
operator was first introduced by G. Aronsson in 60's when he studied
the $L^{\infty}$ variational problems ([A1-4]). We say that $u\in
W_{loc}^{1,\infty}(\Omega)$ is an {\it{absolute minimizer for $H$ in
$\Omega$}} if for any bounded open set $V\subset \bar V\subset
\Omega$ and $v\in W^{1,\infty}(V)$,
$$
u|_{\partial V}=v|_{\partial V}
$$
implies that
$$
\mathrm{esssup_{V}}H(x,u,Du)\leq \mathrm{esssup_{V}}H(x,v,Dv).
$$
Under suitable assumptions on $H$, it was proved that if $u$ is an
absolute minimizer for $H$ in $\Omega$, then it is a viscosity
solution of the Aronsson equation
$$
A_{H}(u)=0  \quad \text{in $\Omega$}.
$$
See for instance Barron-Jensen-Wang [BJW], Crandall [C] and
Crandall-Wang-Yu [CWY]. When $H={1\over 2}|p|^2$, the Aronsson
operator is the famous infinity Laplacian operator $\Delta_\infty
u=u_{x_i}u_{x_j}u_{x_ix_j}$. We refer to the {\it User's Guide},
Crandall-Ishii-Lions [CIL], for definitions of viscosity solutions.

In a recent interesting paper, Peres-Schramm-Sheffield-Wilson
[PSSW], the authors derived the infinity Laplacian operator from a
two-player zero-sum game, called ``tug-of-war". Roughly speaking,
for fixed $\epsilon>0$, starting from $x_0\in \Omega$, at the $k$-th
turn, the players toss a coin and the winner chooses an $x_k\in
\Omega$ with $|x_k-x_{k-1}|\leq \epsilon$. The game ends when
$x_k\in \partial \Omega$. Player I tries to maximize its payoff
$F(x_k)+{\epsilon ^2\over 2}\sum_{i=0}^{k-1}f(x_i)$ and player II
tries to minimize it. Here $F\in C(\partial \Omega)$ is the
{\it{terminal payoff function}} and $f\in C(\Omega)$ is the {\it
running payoff function}. The setting of [PSSW] is in a length
space. In this paper, we only consider $\Rset ^n$. Let
$u_{\epsilon}$ be the value of the
 above game. Under proper assumptions on $\Omega$, $f$ and $F$,
for example if $\Omega$ is bounded and $f$ is positive, it was
proved in [PSSW] that
$$
\lim_{\epsilon\to 0}u_{\epsilon}=u,
$$
where $u$ is the unique viscosity solution of the following equation
\beq{} {\Delta_{\infty}u\over |Du|^2}=-f  \quad \text{in $\Omega$}
\eeq and
$$
u|_{\partial \Omega}=F.
$$
Since $|Du|$ might be zero in equation (1.1), the definition of a
viscosity solution of equation (1.1) is little bit subtle. At the
touching point where the gradient of a test function $\phi$
vanishes, we need to consider $\max_{v\in S^{n-1}}v\cdot
D^2\phi\cdot v$ or $\min_{v\in S^{n-1}}v\cdot D^2\phi\cdot v$
depending on whether $\phi$ touches from above or from below. See
Definition 2.1 in next section. It is natural to multiply both sides
of equation (1.1) by $|Du|^2$ to get another equation which looks
nicer \beq{} \Delta_{\infty}u+f(x)|Du|^2=0 \quad \text{in $\Omega$
}. \eeq Here we want to remark that these two equations are not
equivalent. It is easy to show that any viscosity solution of (1.1)
is also a viscosity solution of equation (1.2). However, except when
$f(x)\equiv 0$, a viscosity solution of equation (1.2) might not be
a viscosity solution of equation (1.1). For example, $u\equiv 0$ is
a smooth solution of $(u')^2u''=(u')^2$, but not a solution of
${(u')^2u''\over (u')^2}=1$.

In Barron-Evans-Jensen [BEJ], the authors considered generalized
``tug-of-war" games where the movement of two players satisfies
other dynamics. They derived that the resulting value functions
satisfy PDEs which involve Aronsson-type operators.  They also
provided several other interesting contexts where Aronsson operators
arise. A basic question about the Aronsson operator is whether it is
single valued. Precisely speaking, assume that $u\in C(\Omega)$ is a
viscosity solution of two equations
$$
A_{H}(u)=f(x),  \quad A_{H}(u)=g(x)  \quad \text{in $\Omega$}.
$$
Do we have that $f=g$? In this paper, we show that the answer is
``Yes" if $H\in C^2(\Omega\times \Rset \times \Rset ^n)$ and $u\in
W_{loc}^{1,\infty}(\Omega)$. The assumption that $|Du|$ is locally
bounded is not necessary for a large class of  $H$. This conclusion
that Aronsson operator has unique value is not obvious at all since
$u$
 lacks sufficient regularity. For example, $u=x^{4\over 3}-y^{4\over 3}$ is a viscosity solution
  of the infinity Laplacian equation, but it is only $C^{1,{1\over 3}}$. We will in fact
  prove a more general result.

  \bth{} Assume that $B\in C^1(\Omega \times \Rset
\times \Rset ^n, \Rset ^n)$, $c\in C(\Omega\times \Rset \times \Rset
^n)$ and $f,g\in C(\Omega)$. Suppose that $u\in
W^{1,\infty}_{loc}(\Omega)$ is a viscosity subsolution of the
following equation \beq{} B(x,u,Du)\cdot D^2u \cdot
B(x,u,Du)+c(x,u,Du)=f(x), \eeq and is a viscosity supersolution of
the following equation \beq{} B(x,u,Du)\cdot D^2u \cdot
B(x,u,Du)+c(x,u,Du)=g(x).\eeq  Then
$$
f\leq g   \quad \text{in $\Omega$}.
$$
\eth

It is clear that if $H\in C^2$, the Aronsson operator $A_{H}$
satisfies structure assumptions in Theorem 1.1. Hence the Aronsson
operator is single valued for locally Lipschitz continuous viscosity
solutions. For suitable $H$, including the most interesting case
$H={1\over 2}|p|^2$, the assumption $u\in
W_{loc}^{1,\infty}(\Omega)$ can be relaxed to $u\in C(\Omega)$. See
Corollary 3.2 and Remark 3.3.

In [PSSW], the authors proposed an open problem which asks whether
two different running payoff functions will lead to the same value
function. See Problem 4 at the end of [PSSW]. We will show that the
answer is No. The following is the precise statement. \bth{} Assume
that $f,g\in C(\Omega)$. Suppose that $u\in C(\Omega)$ is
simultaneously a viscosity solution of two equations \beq{}
{\Delta_{\infty}u\over |Du|^2}=-g(x), \quad {\Delta_{\infty}u\over
|Du|^2}=-f(x) \quad \text{in $\Omega$}. \eeq
 Then
$$
f=g.
$$
\eth

We want to stress that the above theorem can not be deduced directly
from Theorem 1.1. In fact, its proof is much more tricky. We need to
employ the endpoint estimate (2.2) developed in [CEG] to avoid the
situation where $|Du|$ vanishes. If we want to apply Theorem 1.1 to
prove Theorem 1.2, we need an open subset of $\Omega$ where $|Du|$
is bounded away from $0$. The existence of such an open subset
requires the continuity of $|Du|$, which can be given a meaning
independent of the existence of $Du$ itself. The author has proved
this continuity if $n=2$, but has no clue how to prove in higher
dimensions.

\vspace{2mm}

We note that the question of whether or not a function $u$ can
simultaneously solve two distinct Hamilton-Jacobi equations
$H(Du)=f$ and $H(Du)=g$ in the viscosity sense is mentioned in
[CIL]. If $n=1$ or $H$ is uniformly continuous, it was proved in
Evans [E] that the answer is No. Frankowska [F] also provided some
sufficient conditions on $H$ which lead to $f=g$. However for
general situations, this question remains open.

\vspace{2mm}

Our paper is organized as follows. In Section 2, we will review some
known results in [CEG]. In Section 3, we will prove Theorem 1.1 and
Theorem 1.2.

\vspace{2mm}

{\bf Acknowledgement} The author would like to thank Prof. Michael
 Crandall for many valuable comments and suggestions which
significantly  improved our presentation of these results. His
encouragement, as always, is deeply appreciated.

\sectionnew{Preliminaries}
Viscosity solutions of equation (1.1) are defined as follows. \bde{}
$u\in C(\Omega)$ is a viscosity supersolution of equation (1.1) in
$\Omega$ if for any $x_0\in \Omega$ and $\phi\in C^2(\Omega)$
satisfying
$$
0=\phi (x_0)-u(x_0)\geq \phi (x)-u(x)  \quad \text{for all $x\in
\Omega$},
$$
one of the following holds:\\
(1) $D\phi (x_0)\ne 0$ and
$$
\Delta_{\infty} \phi (x_0)\leq -f(x_0)|D\phi (x_0)|^2;
$$
or,\\
(2) $D\phi (x_0)=0$ and
$$
\min_{\{p\in \Rset ^n|\ |p|=1\}}p\cdot D^2\phi (x_0)\cdot p\leq
-f(x_0).
$$
$u\in C(\Omega)$ is a viscosity subsolution of equation (1.1) in
$\Omega$ if for any $x_0\in \Omega$ and $\phi\in C^2(\Omega)$
satisfying
$$
0=\phi (x_0)-u(x_0)\leq \phi (x)-u(x)  \quad \text{for all $x\in
\Omega$},
$$
one of the following holds:\\
(1) $D\phi (x_0)\ne 0$ and
$$
\Delta_{\infty} \phi (x_0)\geq -f(x_0)|D\phi (x_0)|^2;
$$
or,\\
(2) $D\phi (x_0)=0$ and
$$
\max_{\{p\in \Rset ^n|\ |p|=1\}}p\cdot D^2\phi (x_0)\cdot p\geq
-f(x_0).
$$
$u$ is a viscosity solution of equation (1.1) in $\Omega$ if it is
both a viscosity supersolution and subsolution. \ede

A very useful tool to study the infinity Laplacian operator is
``comparison with cones" which was introduced in [CEG](see
definition 2.3). In this terminology, it had been proved in Jensen
[J] that viscosity supersolutions (subsolutions) of the infinity
Laplacian equation enjoy comparison with cones from blow
(respectively, above), and that if $u\in C(\Omega)$ enjoys
comparison with cones from above or from blow, then $u\in
W_{loc}^{1,\infty}(\Omega)$. Crandall, Evans and Gariepy went on to
observe that if $u$ is upper semicontinuous and enjoys comparison
with cones from above, then it is a subsolution of the infinity
Laplacian equation and the quantity
$$
{\max_{\partial B_{r}(x)}u-u(x)\over r}
$$
is nondecreasing with respect to $r$. Hence one can define
$$
S_{u,+}(x)=\lim_{r\to 0^{+}}{\max_{\partial B_{r}(x)}u-u(x)\over r}.
$$
It turns out the function $S_{u,+}(x)$ has the following properties:\\
(1) $S_{u,+}(x)$ is upper-semicontinuous and \beq{}
S_{u,+}(x)=\lim_{r\to 0}\mathrm{esssup}_{B_r(x)}|Du| \eeq
(2) If $u$ is differentiable at $x$, then
$S_{u,+}(x)=|Du(x)|$.\\[2mm]
(3) (Endpoint estimate.) Assume that $x_r\in \partial B_{r}(x)$ and
$u(x_r)=\max_{\partial B_{r}(x)}u$, then \beq{} S_{u,+}(x_r)\geq
{\max_{\partial B_{r}(x)}u-u(x)\over r}\geq S_{u,+}(x).\eeq

\vspace{2mm}

Other notations we have used or will use includes:

\vspace{2mm}

$\bullet$ $\Omega$ is a bounded open subset of $\Rset ^n$.

$\bullet$ For any set $V\in \Rset ^n$, $\partial V$ is its boundary
and $\bar V$ is its closure.

$\bullet$ $B_{r}(x)$ is the open ball $\{y\in \Rset ^n|\ |y-x|<r\}$,
where $|\cdot|$ is the Euclidean norm.

$\bullet$ $S^{n-1}$ denotes the unit sphere in $\Rset ^n$.

$\bullet$ For $p$, $q\in \Rset ^n$, $p\cdot q$ is the usual inner
product of $p$ and $q$. If $A$ is an $n\times n$ matrix, the $p\cdot
A\cdot q$ means $p\cdot (Aq)$.

$\bullet$ For $p\in \Rset ^n$, $p\otimes p$ is the $n\times n$
matrix whose $(i,j)$ entry is $p_ip_j$.

$\bullet$ If $f:\Omega\to \Rset$,  then $Df$ is its gradient and
$D^{2}f$ is  its Hessian matrix.

\sectionnew{Proofs}
To prove Theorem 1.1, we first prove the following lemma. Our proof
heavily depends on the highly degenerate structure of equations
(1.3) and (1.4) and an elegant inequality in [CIL].

\ble{}Let $\tau_1, \tau_2\in\R,$ $\tau_2<\tau_1,$ and the
assumptions on $B, c$ of Theorem 1.1 be satisfied. Then there does
not exist a Lipschitz continuous function $u$ in $\bar \Omega$ such
that the following three conditions hold:
\begin{itemize}
\item[(i)] $u$ is a viscosity subsolution of
\begin{equation*}B(x,u,Du)\cdot D^2u\cdot
B(x,u,Du)+c(x,u,Du)=\tau_1\text{\ in\ } \Omega,
\end{equation*}
\item[(ii)]
$u$ is a viscosity supersolution of
\begin{equation*}B(x,u,Du)\cdot D^2u\cdot
B(x,u,Du)+c(x,u,Du)=\tau_2\text{\ in\ } \Omega,
\end{equation*}
\item[(iii)]\hskip2in
$ u=0 \text{\ on\ }\partial \Omega$.
\end{itemize}  \ele \proof We argue by contradiction. Suppose that  there exists
a Lipschitz continuous function $u$ in $\bar \Omega$ which satisfies
(i)-(iii). Let us denote \beq{} \sup _{x\ne y\in \bar
\Omega}{|u(x)-u(y)|\over |x-y|}=C<+\infty; \eeq Without loss of
generality,
 we assume that there exists some $x_0\in \Omega$ such that $u(x_0)=1$. For $\epsilon>0$, let
 $$
 u_{\epsilon}=(1+\epsilon ^{3\over 4})u
 $$
 and
 $$
 w_{\epsilon}(x,y)=u_{\epsilon}(x)-u(y)- {1\over 2\epsilon}|x-y|^2 .
 $$
 Let $(\bar x, \bar y)\in \bar \Omega\times \bar \Omega$ such that
 $$
 w_{\epsilon}(\bar x, \bar y)=\max_{ \bar \Omega\times \bar \Omega}w_{\epsilon}.
 $$
 Owing to (3.1), we have that $|\bar x-\bar y|=O(\epsilon)$. By  (iii) and (3.1),
 if $(\bar x, \bar y)\in \partial  (\Omega\times \Omega)$,
 then $w_{\epsilon}(\bar x, \bar y)=O(\epsilon)$. Note that  $w_{\epsilon}(x_0,x_0)=\epsilon ^{3\over 4}$.
 Hence when $\epsilon$ is small enough, $(\bar x, \bar y)\in \Omega\times \Omega$.
 According to Crandall-Ishii-Lions [CIL], there exist two $n \times n$ symmetric matrices $X$ and $Y$ such that
$$
({1\over {\epsilon}}(\hat x-\hat y),X)\in {\bar J}_{V}^{2,+}u_{\epsilon}(\hat
x),\ ({1\over {\epsilon}}(\hat x-\hat y),Y)\in {\bar
J}_{V}^{2,-}u(\hat y)
$$
and
\beq{}
-{3\over \epsilon}
\begin{pmatrix}
I_{n} & 0\\
0 & I_n
\end{pmatrix}
\leq
\begin{pmatrix}
X & 0\\
0 & -Y
\end{pmatrix}
\leq {3\over \epsilon}
\begin{pmatrix}
I_{n} & -I_{n}\\
-I_{n} & I_n
\end{pmatrix}
, \eeq See [CIL] for definitions of ${\bar J}_{V}^{2,+}$ and ${\bar
J}_{V}^{2,-}$. It is easy
to see that $u_{\epsilon}$ is a viscosity subsolution of
$$
B(x,{u_{\epsilon}\over {1+{\epsilon ^{3\over 4}}}},
{Du_{\epsilon}\over {1+{\epsilon ^{3\over 4}}}})\cdot
D^2u_{\epsilon}\cdot B(x,{u_{\epsilon}\over {1+{\epsilon ^{3\over
4}}}},  {Du_{\epsilon}\over {1+{\epsilon ^{3\over
4}}}})+{{(1+{\epsilon ^{3\over 4}})}} c(x,{u_{\epsilon}\over {1+{\epsilon ^{3\over 4}}}},
{Du_{\epsilon}\over {1+{\epsilon ^{3\over 4}}}})={\tau_{1}
{(1+{\epsilon ^{3\over 4}})}}.
$$
Hence \beq{} B\mathrm {\Big (}\bar x,u(\bar x),  {\bar x-\bar y\over
{\epsilon}{(1+{\epsilon ^{3\over 4}})}}\mathrm {\Big )}\cdot X\cdot
B\mathrm {\Big (}\bar x, u(\bar x),  {\bar x-\bar y\over \epsilon
{(1+{\epsilon ^{3\over 4}})}}\mathrm {\Big )}+{{(1+{\epsilon
^{3\over 4}})}} c\mathrm {\Big (}\bar x,u(\bar x), {\bar x-\bar
y\over \epsilon {(1+{\epsilon ^{3\over 4}})}}\mathrm {\Big )}\geq
{\tau_{1} {(1+{\epsilon ^{3\over 4}})}}. \eeq By (ii), \beq{}
B\mathrm {\Big (}\bar y,u(\bar y), {\bar x-\bar y\over
{\epsilon}}\mathrm {\Big )}\cdot Y\cdot B\mathrm {\Big (}\bar
y,u(\bar y), {\bar x-\bar y\over {\epsilon}}\mathrm {\Big
)}+c\mathrm {\Big (}\bar y,u(\bar y),{\bar x-\bar y\over
{\epsilon}}\mathrm {\Big )}\leq \tau_{2}. \eeq Owing to the right
hand side inequality in (3.2), we have that for $v_1,v_2\in \Rset
^n$,
$$
v_1\cdot X\cdot v_1-v_2\cdot Y\cdot v_2\leq {3\over {\epsilon}}|v_1-v_2|^2.
$$
Choosing
$$
v_1=B\mathrm {\Big (}\bar x,u(\bar x),  {\bar x-\bar y\over
{\epsilon}{(1+{\epsilon ^{3/4}})}}\mathrm {\Big )},\ v_2=B\mathrm
{\Big (}\bar y,u(\bar y), {\bar x-\bar y\over {\epsilon}}\mathrm
{\Big )}
$$
and using $|\bar x-\bar y|=O(\epsilon)$, (3.3)-(3.4), we discover
that
$$
o(1)\geq {\tau_1 {(1+{\epsilon ^{3\over
4}})}}-\tau_2,
$$
where $\lim_{\epsilon \to 0}o(1)=0$. This is impossible when
$\epsilon$ is small enough. So Lemma 3.1 holds.\qed

\vspace{2mm}

{\bf{Proof of Theorem 1.1}}. For any $x_0\in \Omega$, if
$g(x_0)<f(x_0)$, then there exists $r>0$ and $\tau_1>\tau_2$ such
that $\bar B_r(x_0)\subset \Omega$ and
$$
g(x)<\tau_2<\tau_1<f(x)  \quad \text{in $B_r(x_0)$}.
$$
Choose $K$ large enough such that
$$
u(x)< u(x_0)+Kr^2  \quad \text{on $\partial B_r(x_0)$}.
$$
Denote $\delta={1\over 2}\min_{\partial B_r(x_0)}(u(x_0)+Kr^2-u(x))$
and $v(x)=u(x)-u(x_0)-K|x-x_0|^2+\delta$. Let
$$
V=\{x\in B_{r}(x_0)|\ v(x)> 0\}.
$$
Obviously, $\bar V\subset B_{r}(x_0)$. If we define
$$
\tilde B(x,z,p)=B(x,z+u(x_0)+K|x-x_0|^2-\delta,p+2K(x-x_0))
$$
and
$$
\tilde c(x,z,p)=c(x,z+u(x_0)+K|x-x_0|^2-\delta,p+2K(x-x_0))+2K|\tilde
B(x,z,p)|^2,
$$
$v$ is a viscosity subsolution of
$$
\tilde B(x,v,Dv)\cdot D^2v\cdot \tilde B(x,v,Dv)+\tilde c(x,v,Dv)=\tau_1   \quad \text{in $V$}
$$
and a viscosity supersolution of
$$
\tilde B(x,v,Dv)\cdot D^2v\cdot \tilde B(x,v,Dv)+\tilde c(x,v,Dv)=\tau_2    \quad \text{in $V$}.
$$
Note that in the open set $V$, $v$ satisfies (i)-(iii) in Lemma 3.1.
Since $\bar V\subset \Omega$, $u$ is Lipschitz continuous in $\bar
V$. Hence $v$ is also Lipschitz continuous in $\bar V$. This is a
contradiction. Hence $g(x_0)\geq f(x_0)$. So $f\leq g$.\qed

\bco{} Suppose that $u,f,g\in C(\Omega)$ and $u$ is simultaneously a
viscosity solution of two equations \beq{} {\Delta_{\infty}u}=f(x),
\quad {\Delta_{\infty}u}=g(x) \quad \text{in $\Omega$}. \eeq Then
$$
f=g.
$$
\eco

\proof We argue by contradiction. If not, then there exists $x_0\in
\Omega$ such that $f(x_0)\ne g(x_0)$. Without loss of generality, we
may assume that $f(x_0)>g(x_0)$. Then one of the following must
occur: (i) $f(x_0)>0$, (ii) $g(x_0)<0$. Let us first look at case
(i). Since Corollary 3.2 is a local problem, we may assume that
$$
f(x)>\max\{0,g(x)\}  \quad \text{for $x\in \Omega$}.
$$
 Hence $u$ is a viscosity subsolution of the infinity Laplacian equation
$$
\Delta_{\infty}u=0         \quad \text{in $\Omega$}.
$$
According to [CEG], $u\in W^{1,\infty}_{loc}(\Omega)$. Hence by
Theorem 1.1, $f=g$ in $\Omega$. This is a contradiction. Hence case
(i) will not occur. Similarly, we can show that case (ii) will not
occur either. This is a contradiction. Hence the above corollary
holds.\qed

\bre{} Gariepy-Wang-Yu [GWY], Yu [Y] and Juutinen [Ju] provided a
class of Hamiltonians $H$ such that if $u\in C(\Omega)$ is a
viscosity subsolution or a viscosity supersolution of the Aronsson
equation
$$
A_{H}(u)=0  \quad \text{in $\Omega$},
$$
then  $u\in W_{loc}^{1,\infty}(\Omega)$. Hence by the proof of
Corollary 3.2,  for those $H$, the Aronsson operator $A_{H}$ is also
single valued under the weaker assumption $u\in C(\Omega)$.
 \ere \qed

To prove Theorem 1.2, we first prove the following lemma.

\ble{} Suppose that $\tau_1\ne \tau_2$. Then $u\in W^{1,\infty}(
B_1(0))$ can not be simultaneously a viscosity solution of two
equations \beq{} {\Delta_{\infty}u\over |Du|^2}=\tau_1, \quad
{\Delta_{\infty}u\over |Du|^2}=\tau_2.\eeq \ele

\proof We argue by contradiction. Without loss of generality, let us
assume that $\tau_1>\max\{0,\tau_2\}$. According to the definition
of solutions of
$$
{\Delta_{\infty}u\over |Du|^2}=\tau_1>0,
$$
$u$ can not be constant in any open subset of $B_1(0)$. So we may
assume that $|Du|(0)=\delta>0$, where $|Du|(x)=S_{u,+}(x)$. Let us
denote \beq{} \mathrm{esssup}_{B_1(0)}|Du|=C.\eeq
 Consider
$$
w_{\epsilon}(h)=\max_{x,y\in \bar B_{1\over
2}(0)}(u(x+h)-u(y)-|y|^4-{1\over 2\epsilon}|x-y|^2).
$$
 Choose $h_{\epsilon}\in \partial
B_{{\epsilon^{3/4}}}(0)$ such that
$$
w_{\epsilon}(h_{\epsilon})=\max_{h\in \partial B_{{\epsilon^{3/
4}}}(0)}w_{\epsilon}.
$$
Suppose that
$$
w_{\epsilon}(h_{\epsilon})=u(x_{\epsilon}+h_{\epsilon})-u(y_{\epsilon})-|y_{\epsilon}|^4-{1\over
2\epsilon}|x_{\epsilon}-y_{\epsilon}|^2
$$
for $x_{\epsilon}$, $y_{\epsilon}\in \bar B_{1\over 2}$. Owing to
(3.7), it is not hard to show that when $\epsilon$ is small \beq{}
{1\over \epsilon}|x_{\epsilon}-y_{\epsilon}|\leq C, \quad
|y_{\epsilon}|^4\leq K{\epsilon}^{3\over 4}, \eeq where $K$ is a
constant independent of $\epsilon$. Moreover, it is clear that\beq{}
u(x_{\epsilon}+h_{\epsilon})=\max_{h\in
\partial B_{{\epsilon^{3/4}}}(0)}u(x_{\epsilon}+h)=\max_{y\in \partial B_{
{\epsilon^{3/4}}}(x_{\epsilon})}u(y).\eeq Owing to the definition of
$w_{\epsilon}(h)$, \beq{}
u(x_{\epsilon}+h_{\epsilon})-u(y_{\epsilon})-|y_{\epsilon}|^4-{1\over
2\epsilon}|x_{\epsilon}-y_{\epsilon}|^2\geq \max_{h\in B_{
{\epsilon^{3/4}}}(0)}(u(h)-u(0)). \eeq Also,
\begin{eqnarray*}
{u(x_{\epsilon}+h_{\epsilon})-u(x_{\epsilon})\over {\epsilon^{3\over
4}}}&{\displaystyle={u(x_{\epsilon}+h_{\epsilon})-u(y_{\epsilon})\over
{\epsilon^{3\over 4}}}+{u(y_{\epsilon})-u(x_{\epsilon})\over
{\epsilon^{3\over 4}}}}\\[5mm]
&{\displaystyle\geq
{u(x_{\epsilon}+h_{\epsilon})-u(y_{\epsilon})\over {\epsilon^{3\over
4}}}-{C|x_{\epsilon}-y_{\epsilon}|\over
{\epsilon^{3\over 4}}}}\\[5mm]
&{\displaystyle\geq
{u(x_{\epsilon}+h_{\epsilon})-u(y_{\epsilon})\over {\epsilon^{3\over
4}}}-C^2{\epsilon^{1\over 4}}}.
\end{eqnarray*}
According to (3.10),
$$
{u(x_{\epsilon}+h_{\epsilon})-u(y_{\epsilon})\over {\epsilon^{3\over
4}}}\geq \max_{h\in \partial B_{{\epsilon^{3/
4}}}(0)}{u(h)-u(0)\over {\epsilon^{3\over 4}}}\geq |Du|(0)=\delta.
$$
Since $u$ is a viscosity subsolution of the infinity Laplacian
equation, due to the endpoint estimate (2.2), \beq{}
|Du|(x_{\epsilon}+h_{\epsilon})\geq
{u(x_{\epsilon}+h_{\epsilon})-u(x_{\epsilon})\over {\epsilon^{3\over
4}}}\geq \delta-C^2{\epsilon^{1\over 4}}.\eeq Therefore, when
$\epsilon$ is small, \beq{} |Du|(x_{\epsilon}+h_{\epsilon})\geq
{1\over 2}\delta.\eeq Obviously, when $\epsilon$ is small, both
$x_{\epsilon}$ and $y_{\epsilon}$ are in the interior of $B_{1\over
2}(0)$. According to the {\it User's Guide} Crandall-Ishii-Lions
[CIL], there exist two $n \times n$ symmetric matrices $X$ and $Y$
such that
$$
({1\over {\epsilon}}(x_{\epsilon}-y_{\epsilon}),X)\in {\bar
J}_{V}^{2,+}u(x_{\epsilon}+h_{\epsilon}),\ ({1\over {\epsilon}}(
x_{\epsilon}-y_{\epsilon}),Y)\in {\bar
J}_{V}^{2,-}(u(y_{\epsilon})+|y_{\epsilon}|^4)
$$
and \beq{} -{3\over \epsilon}
\begin{pmatrix}
I_{n} & 0\\
0 & I_n
\end{pmatrix}
\leq
\begin{pmatrix}
X & 0\\
0 & -Y
\end{pmatrix}
\leq {3\over \epsilon}
\begin{pmatrix}
I_{n} & -I_{n}\\
-I_{n} & I_n
\end{pmatrix}
, \eeq See [CIL] for definitions of ${\bar J}_{V}^{2,+}$ and ${\bar
J}_{V}^{2,-}$. Owing to the definition of $|Du|(x)=S_{u,+}(x)$, it
is clear that \beq{} {1\over
{\epsilon}}|x_{\epsilon}-y_{\epsilon}|\geq
|Du|(x_{\epsilon}+h_{\epsilon})\geq {\delta\over 2}. \eeq Since
$u(\cdot+h_{\epsilon})$ is a viscosity solution of equation (3.6),
we have that \beq{} {1\over
{\epsilon}}(x_{\epsilon}-y_{\epsilon})\cdot X\cdot {1\over
{\epsilon}}(x_{\epsilon}-y_{\epsilon})\geq \tau_1{1\over
{\epsilon^2}}|x_{\epsilon}-y_{\epsilon}|^2. \eeq Also,
$$
({1\over
{\epsilon}}(x_{\epsilon}-y_{\epsilon})-4|y_{\epsilon}|^2y_{\epsilon},
Y-4|y_{\epsilon}|^2I_n-8y_{\epsilon}\otimes y_{\epsilon})\in {\bar
J}_{V}^{2,-}u(y_{\epsilon}).
$$
Due to equation (3.6),
$$
({1\over
{\epsilon}}(x_{\epsilon}-y_{\epsilon})-4|y_{\epsilon}|^2y_{\epsilon})\cdot
(Y-4|y_{\epsilon}|^2I_n-8y_{\epsilon}\otimes y_{\epsilon})\cdot
({1\over
{\epsilon}}(x_{\epsilon}-y_{\epsilon})-4|y_{\epsilon}|^2y_{\epsilon})\leq
\tau_2|{1\over
{\epsilon}}(x_{\epsilon}-y_{\epsilon})-4|y_{\epsilon}|^2y_{\epsilon}|^2.
$$
Hence owing to (3.8), \beq{} ({1\over
{\epsilon}}(x_{\epsilon}-y_{\epsilon})-4|y_{\epsilon}|^2y_{\epsilon})\cdot
Y\cdot ({1\over
{\epsilon}}(x_{\epsilon}-y_{\epsilon})-4|y_{\epsilon}|^2y_{\epsilon})\leq
\tau_2{1\over {\epsilon^2}}|x_{\epsilon}-y_{\epsilon}|^2+o(1), \eeq
where $\lim_{\epsilon\to 0}o(1)=0$. Owing to the right hand side
inequality in (8), we have that for $v_1,v_2\in \Rset ^n$,
$$
v_1\cdot X\cdot v_1-v_2\cdot Y\cdot v_2\leq {3\over
{\epsilon}}|v_1-v_2|^2.
$$
Choosing
$$
v_1={1\over {\epsilon}}(x_{\epsilon}-y_{\epsilon}), \ v_2={1\over
{\epsilon}}(x_{\epsilon}-y_{\epsilon})-4|y_{\epsilon}|^2y_{\epsilon}
$$
and using (3.14)-(3.16), one finds that \beq{}
{48|y_{\epsilon}|^6\over {\epsilon}}\geq (\tau_1-\tau_2){1\over
{\epsilon ^2}}|x_{\epsilon}-y_{\epsilon}|^2-o(1)\geq
(\tau_1-\tau_2)({\delta \over 2})^2-o(1). \eeq Owing to (3.8),
$$
{4|y_{\epsilon}|^6\over {\epsilon}}\leq {\epsilon^{1\over 8}}.
$$
This contradicts to (3.17) when $\epsilon$ is small. \qed

\vspace{2mm}

 {\bf{Proof of Theorem
1.2}}. We argue by contradiction. If not, then there exists $x_0\in
\Omega$ such that $f(x_0)\ne g(x_0)$. Without loss of generality, we
may assume that $f(x_0)>g(x_0)$. Then one of the following must
occur: (i) $f(x_0)>0$, (ii) $g(x_0)<0$. Let us first look at case
(i). Since Theorem 1.2 is a local result, we may assume that
$$
f(x)>\tau_1>\tau_2>\max\{0,g(x)\}  \quad \text{for $x\in \Omega$},
$$
where $\tau_1$ and $\tau_2$ are two positive constants. Hence $u$ is
also a viscosity supersolution of the infinity Laplacian equation
$$
\Delta_{\infty}u=0 \quad \text{in $\Omega$}.
$$
According to [CEG], $u\in W_{loc}^{1,\infty}(\Omega)$. Choose $r>0$
such that $\overline { B_r(x_0)}\subset \Omega$. Then $u\in
W^{1,\infty}(B_r(x_0))$. Consider
$$
u_{r}(x)=u(rx+x_0).
$$
Then $u_r\in W^{1,\infty}(B_1(0))$ and it is simultaneously a
viscosity solution of two equations\beq {\Delta_{\infty}u_r\over
{|Du_r|^2}}=-r^2\tau_1, \quad {\Delta_{\infty}u_r\over
{|Du_r|^2}}=-r^2\tau_2 \quad \text{in $B_1(0)$}\eeq  This
contradicts to Lemma 3.4. Similarly, we can show that case (ii) will
not happen either. Hence Theorem 1.2 holds.\qed

\bre{} By obvious modifications, our method can be used to prove
that any operator like $\Delta_{\infty}u/f(x,u,Du)$ is single valued
if $f\in C(\Omega\times \Rset \times \Rset ^n)$ is nonnegative or
nonpositive. Here is a potential application.  A very interesting
problem about the infinity Laplacian operator is to find its
geometric interpretation. More specifically, does there exist a
function $f$ such that $\Delta_{\infty}u/f(x,u,Du)$ represents some
kind of curvature of the graph of $u$ (might be in viscosity sense)?
The answer of this question will justify the study of the parabolic
infinity Laplacian equation. Our results implies that if there is
indeed such a curvature, then it is well defined, i.e, a surface has
at most one curvature. \ere

\renewcommand{\em}{\textrm}
\begin{small}
  \renewcommand{\refname}{ {\flushleft\normalsize\bf{References}} }

\end{small}

\end{document}